\documentstyle[12pt]{amsart}
\def\opn#1#2{\def#1{\operatorname{#2}}} 
\opn\chara{char}
\opn\length{\ell}
\opn\pd{pd}
\opn\rk{rk}
\opn\projdim{proj\,dim}
\opn\injdim{inj\,dim}
\opn\rank{rank}
\opn\depth{depth}
\opn\grade{grade}
\opn\height{height}
\opn\embdim{emb\,dim}
\opn\codim{codim}

\opn\Tr{Tr}
\opn\bigrank{big\,rank}
\opn\superheight{superheight}\opn\lcm{lcm}
\opn\trdeg{tr\,deg}%
\opn\reg{reg}
\opn\lreg{lreg}
\opn\skel{skel}
\opn\com{com}
\opn\div{div}
\opn\Div{Div}
\opn\cl{cl}
\opn\Cl{Cl}
\opn\Spec{Spec}
\opn\Supp{Supp}
\opn\supp{supp}
\opn\Sing{Sing}
\opn\Ass{Ass}
\opn\Ann{Ann}
\opn\Rad{Rad}
\opn\Soc{Soc}
\opn\Ker{Ker}
\opn\Coker{Coker}
\opn\Im{Im}
\opn\Hom{Hom}
\opn\Tor{Tor}
\opn\Ext{Ext}
\opn\End{End}
\opn\Aut{Aut}
\opn\id{id}

\opn\nat{nat}
\opn\pff{pf}
\opn\Pf{Pf}
\opn\GL{GL}
\opn\SL{SL}
\opn\mod{mod}
\opn\ord{ord}
\opn\aff{aff}
\opn\con{conv}
\opn\relint{relint}
\opn\st{st}
\opn\lk{lk}
\opn\cn{cn}
\opn\core{core}
\opn\vol{vol}
\opn\link{link}
\opn\star{star}
\opn\gr{gr}

\def\pot#1#2{#1[\kern-0.28ex[#2]\kern-0.28ex]}

\opn\dirlim{\underrightarrow{\lim}}
\opn\inivlim{\underleftarrow{\lim}}
\let\to=\rightarrow

\def\Implies{\ifmmode\Longrightarrow \else
     \unskip${}\Longrightarrow{}$\ignorespaces\fi}
\def\implies{\ifmmode\Rightarrow \else
     \unskip${}\Rightarrow{}$\ignorespaces\fi}
\def\iff{\ifmmode\Longleftrightarrow \else
     \unskip${}\Longleftrightarrow{}$\ignorespaces\fi}

\let\:=\colon
\newtheorem{Theorem}{Theorem}[section]
\newtheorem{Lemma}[Theorem]{Lemma}
\newtheorem{Corollary}[Theorem]{Corollary}

\let\epsilon\varepsilon
\let\phi=\varphi
\let\kappa=\varkappa
\def\qed{\ifhmode\textqed\fi
   \ifmmode\ifinner\quad\qedsymbol\else\dispqed\fi\fi}
\def\textqed{\unskip\nobreak\penalty50
    \hskip2em\hbox{}\nobreak\hfil\qedsymbol
    \parfillskip=0pt \finalhyphendemerits=0}
\def\dispqed{\rlap{\qquad\qedsymbol}}

\opn\initial{in}
\opn\inim{inm}
\opn\rev{rev}
\opn\Gin{Gin}
\opn\Shift{Shift}
\opn\shift{shift}
\opn\rate{rate}
\opn\Mon{Mon}
\opn\lex{lex}
\opn\rev{rev}
\opn\red{red}
\opn\max{max}
\opn\min{min}
\opn\initial{in}
\opn\Ker{Ker}
\opn\GL{GL}
\opn\proj{proj}
\begin{document}
\title{Algebraic shifting and graded Betti numbers}
\author{Satoshi Murai and Takayuki Hibi}
\date{}
\maketitle
\begin{abstract}
Let $S = K[x_1, \ldots, x_n]$ denote
the polynomial ring in $n$ variables over a field 
$K$ with each $\deg x_i = 1$.
Let $\Delta$ be a simplicial complex on 
$[n] = \{ 1, \ldots, n \}$
and $I_\Delta \subset S$ its Stanley--Reisner ideal.  
We write $\Delta^e$ for the exterior algebraic
shifted complex of $\Delta$ and 
$\Delta^c$ for a combinatorial shifted complex 
of $\Delta$.  Let 
$\beta_{ii+j}(I_{\Delta}) = \dim_K \Tor_i(K, I_\Delta)_{i+j}$ 
denote the graded Betti numbers of $I_\Delta$.
In the present paper it will be proved that
(i)  
$\beta_{ii+j}(I_{\Delta^e}) \leq
\beta_{ii+j}(I_{\Delta^c})$
for all $i$ and $j$, 
where the base field is infinite, and
(ii) $\beta_{ii+j}(I_{\Delta}) \leq
\beta_{ii+j}(I_{\Delta^c})$
for all $i$ and $j$,
where the base field is arbitrary.
Thus in particular 
one has
$\beta_{ii+j}(I_\Delta) \leq \beta_{ii+j}(I_{\Delta^{lex}})$
for all $i$ and $j$, where
$\Delta^{lex}$ is the unique lexsegment simplicial complex
with the same $f$-vector as $\Delta$ and where the base field 
is arbitrary.
\end{abstract}

\section*{Introduction}
Kalai \cite{Kalai} together with Herzog \cite{Herzog}
offer an attractive introduction, 
which includes several unsolved problems and conjectures,
to the combinatorial and 
algebraic study of shifting theory in algebraic and 
extremal combinatorics.
 
Let $S = K[x_1, \ldots, x_n]$ denote
the polynomial ring in $n$ variables over 
a field $K$ with each $\deg x_i = 1$.  
One of the current trends 
in computational commutative algebra
is the computation of the graded Betti numbers of
homogeneous ideals.  
Recall that the graded Betti numbers
$\beta_{ij} = \beta_{ij}(I)$, where $i, j \geq 0$, of 
a homogeneous ideal $I \subset S$ are
\[
\beta_{ij}(I) = \dim_K \Tor_i(K, I)_{j}.
\]
In other words, the graded Betti numbers 
$\{ \beta_{ij} \}_{i,j = 0, 1, \ldots}$
appear in the 
minimal graded free resolution 
\[
0 \longrightarrow \bigoplus_{j} S(-j)^{\beta_{hj}}
\longrightarrow
\cdots
\longrightarrow
\bigoplus_{j} S(-j)^{\beta_{1j}}
\longrightarrow
\bigoplus_{j} S(-j)^{\beta_{0j}}
\longrightarrow I
\longrightarrow 0. 
\]
of $I$ over $S$, 
where $h = \proj \dim_S I$ is the projective
dimension of $I$ over $S$.

Let $\Delta$ be a simplicial complex on 
$[n] = \{ 1, \ldots, n \}$ and $I_\Delta \subset S$
the Stanley--Reisner ideal of $\Delta$.  
We write $\Delta^s$,
$\Delta^e$ and $\Delta^c$ for the symmetric algebraic
shifted complex, the exterior algebraic shifted complex
and a combinatorial shifted complex, respectively,
of $\Delta$.
Since the paper \cite{AraHerHiLexsegment}
was published, it has been conjectured that for
an arbitrary simplicial 
complex $\Delta$ on $[n]$ one has
\[
\beta_{ii+j}(I_\Delta) \leq
\beta_{ii+j}(I_{\Delta^s}) \leq
\beta_{ii+j}(I_{\Delta^e}) \leq
\beta_{ii+j}(I_{\Delta^c})
\]
for all $i$ and $j$.
When the base field is of characteristic $0$,
the first inequality
$\beta_{ii+j}(I_\Delta) \leq
\beta_{ii+j}(I_{\Delta^s})$
is proved in \cite[Theorem 2.1]{AraHerHiJAC}.

Let $\Delta'$ be a shifted (or strongly stable
\cite[p.\ 365]{AraHerHiLexsegment}) simplicial complex
with the same $f$-vector as $\Delta$ and $\Delta^{\lex}$ 
the unique lexsegment simplicial complex
with the same $f$-vector as $\Delta$
(\cite[Theorem 3.5]{AraHerHiLexsegment}).
It is known \cite[Theorem 4.4]{AraHerHiLexsegment} that
$\beta_{ii+j}(I_{\Delta'}) \leq
\beta_{ii+j}(I_{\Delta^{\lex}})$.
Since $\Delta^s$ is shifted
with the same $f$-vector as $\Delta$,
when the base field is of characteristic $0$,
one has 
$\beta_{ii+j}(I_\Delta) \leq
\beta_{ii+j}(I_{\Delta^{\lex}})$
for all $i$ and $j$
(\cite[Theorem 2.9]{AraHerHiJAC}).

The main purpose of the present paper is 
to establish two fundamental results stated below 
concerning the graded Betti numbers of 
$I_\Delta$, $I_{\Delta^e}$ and $I_{\Delta^c}$.  

\bigskip

\noindent
{\bf Theorem 2.10.}
{\em
Let the base field be infinite. 
Let $\Delta$ be a simplicial complex,
$\Delta^e$ the exterior algebraic shifted complex of
$\Delta$ and
$\Delta^c$ a combinatorial shifted complex of
$\Delta$.
Then 
\[
\beta_{ii+j}(I_{\Delta^e}) \leq \beta_{ii+j}(I_{\Delta^c})
\] 
for all $i$ and $j$. 
}

\bigskip

\noindent
{\bf Theorem 3.4.}
{\em 
Let the base field be arbitrary.
Let $\Delta$ be a simplicial complex 
and $\Delta^c$ a combinatorial shifted complex
of $\Delta$.  Then
\[
\beta_{i i+j}(I_\Delta) \leq \beta_{i i+j}(I_{\Delta^c})
\]
for all $i$ and $j$.
}

\bigskip

Since $\Delta^c$ is shifted 
with the same $f$-vector as $\Delta$,
it follows from Theorem 3.4 together with
\cite[Theorem 4.4]{AraHerHiLexsegment} that

\bigskip

\noindent
{\bf Corollary 3.5.}
{\em 
Let the base field be arbitrary.
Let $\Delta$ be a simplicial complex 
and $\Delta^{\lex}$ 
the unique lexsegment simplicial complex
with the same $f$-vector as $\Delta$.
Then
\label{corollary}
\[
\beta_{ii+j}(I_\Delta) \leq
\beta_{ii+j}(I_{\Delta^{\lex}})
\]
for all $i$ and $j$.
}

\bigskip

The present paper will be organized as follows.
First of all, following \cite{Herzog}
the fundamental materials on algebraic shifting
will be summarized in Section $1$. 
Second, our proof of Theorem \ref{maintheorem}
will be achieved in Section $2$.
On the other hand,
based on Hochster's formula 
\cite[Theorem 5.5.1]{BrunsHerzog}
to compute graded Betti numbers
of Stanley--Reisner ideals, we will prove 
Theorem \ref{theorem} in Section $3$.

Finally, in Section $4$ 
the bad behavior of graded Betti numbers
of $I_{\Delta^c}$ will be studied.
More precisely, 
since a combinatorial shifted complex of $\Delta$
is not unique, it is natural to ask,
given a simplicial complex $\Delta$, if
there exist combinatorial shifted complexes
$\Delta_{\flat}^c$ and $\Delta_{\sharp}^c$ of $\Delta$ 
such that, for each combinatorial shifted complex
$\Delta^c$ of $\Delta$ and for all $i$ and $j$, one has
\[
\beta_{ii+j}(I_{\Delta_{\flat}^c}) 
\leq \beta_{ii+j}(I_{\Delta^c})
\leq \beta_{ii+j}(I_{\Delta_{\sharp}^c}).
\]
Unfortunately, in general,
the existence of such the combinatorial shifted complexes
$\Delta_{\flat}^c$ and $\Delta_{\sharp}^c$
cannot be expected
(Theorem \ref{nonexistence}).  
Especially, we construct
a simplicial complex $\Delta$ for which
there is no combinatorial shifted complex
$\Delta^c$ of $\Delta$ with
$\Delta^e = \Delta^c$ 
(Corollary
\ref{exteriorcombinatorial}).

\section{Algebraic shifting}
Let $[n] = \{ 1, \ldots, n \}$ and write ${[n] \choose i}$ 
for the set of $i$-element subsets of $[n]$.
Let 
$S = K[x_1, \ldots, x_n]$ denote 
the polynomial ring in $n$ variables over 
a field $K$ with
each $\deg x_i = 1$.  Let $V$ be a vector space over $K$
of dimension $n$ with basis $e_1, \ldots, e_n$ and
$E = \bigoplus_{d=0}^{n} \bigwedge^d (V)$ 
the exterior algebra of $V$.
If $\sigma = \{ j_1, \ldots, j_d \} \in {[n] \choose d}$
with $j_1 < \cdots < j_d$, then
$x_\sigma = x_{j_1} 
\cdots x_{j_d}$
is a squarefree monomial of $S$ of degree $d$
and
$e_\sigma = e_{j_1} 
\wedge \cdots \wedge e_{j_d}
\in \bigwedge^d (V)$ will be called a {\em monomial} of $E$ 
of degree $d$.  

Let $\Delta$ be a simplicial complex on $[n]$.  Thus $\Delta$
is a collection of subsets of $[n]$ such that 
(i) $\{ j \} \in \Delta$ for all $j \in [n]$ and
(ii) if $\tau \subset [n]$ and $\sigma \in \Delta$ 
with $\tau \subset \sigma$, then
$\tau \in \Delta$.  
A {\em face} of $\Delta$ is an element $\sigma \in \Delta$.
The $f$-{\em vector} of $\Delta$
is the vector $f(\Delta) = (f_0, f_1, \ldots)$, where $f_i$ is
the number of faces $\sigma \in \Delta$ with
$|\sigma| = i + 1$.
(For a finite set $\sigma$ the notation $|\sigma|$ stands for 
its cardinality.)
The {\em Stanley--Reisner ideal} of $\Delta$
is the ideal $I_\Delta$ of $S$ generated by 
those squarefree monomials $x_\sigma$ with $\sigma \not\in \Delta$.
The {\em exterior face ideal} of $\Delta$ is the ideal
$J_\Delta$ of $E$ generated by those monomials $e_\sigma$ 
with $\sigma \not\in \Delta$.  

If $I \subset S$ is a squarefree ideal, i.e.,
an ideal generated by squarefree monomials,
with each $x_i \not\in I$, 
then there is a unique simplicial complex $\Delta$ on $[n]$
with $I = I_\Delta$.  
If $I \subset E$ is a monomial ideal, i.e.,
an ideal generated by  monomials,
with each $e_i \not\in I$, 
then there is a unique simplicial complex $\Delta$ on $[n]$
with $I = J_\Delta$. 

A monomial ideal $I \subset S$ is called {\em strongly stable}
if for each monomial $u \in I$ and for each $j \in [n]$
for which $x_j$ divides $u$ one has
$x_i u / x_j \in I$ for all $i < j$.
A squarefree ideal $I \subset S$ is called 
{\em squarefree strongly stable}
if for each monomial $x_\sigma \in I$ and for each $j \in \sigma$
one has 
$x_{
(\sigma \setminus \{ j \}) \bigcup \{ i \} 
} \in I$ for all $i < j$ with
$i \not\in \sigma$.
A monomial ideal $I \subset E$ is called {\em strongly stable}
if for each monomial $e_\sigma \in I$ and for each $j \in \sigma$
one has 
$e_{
(\sigma \setminus \{ j \}) \bigcup \{ i \} 
} \in I$ for all $i < j$ with
$i \not\in \sigma$.

We say that a simplicial complex $\Delta$ on $[n]$ 
is {\em shifted}
if the monomial ideal $J_\Delta$ is strongly stable
(or equivalently, the squarefree ideal $I_\Delta$ is
squarefree strongly stable). 
In other word, 
$\Delta$ is shifted if $\Delta$ possesses the property that
for each face $\sigma \in \Delta$ and for each $i \in \sigma$
one has $(\sigma \setminus \{ i \}) \bigcup \{ j \} \in \Delta$
for all $j > i$ with $j \not\in \sigma$. 

Assume that the base field $K$ is of characteristic $0$.
Fix the reverse lexicographic order $<_{\rev}$ 
on $S = K[x_1, \ldots, x_n]$ 
induced by the ordering $x_1 > \cdots > x_n$.
Given a homogeneous ideal $I \subset S$, we write $\Gin^S(I)$
for the {\em generic initial ideal} 
\cite[p.\ 129]{Green}
of $I$ with respect to $<_{\rev}$.
The generic initial ideal $\Gin^S(I)$ of a homogeneous ideal
$I \subset S$ is strongly stable
\cite[Theorem 1.27]{Green}.

We refer the reader to \cite{AraHerHiGotzmann} for the foundation 
on the Gr\"obner basis theory in the exterior algebra. 
Assume that the base field $K$ is infinite.
We work with the reverse lexicographic order $<_{\rev}$
on $E$ induced by the ordering $e_1 > e_2 > \cdots > e_n$.
Given a homogeneous ideal $I \subset E$, we write $\Gin^E(I)$
for the {\em generic initial ideal} 
\cite[p.\ 183]{AraHerHiGotzmann}
of $I$ with respect to $<_{\rev}$.
The generic initial ideal $\Gin^E(I)$ of a homogeneous ideal
$I \subset E$ is strongly stable
\cite[Proposition 1.7]{AraHerHiGotzmann}.

A {\em shifting operation} on $[n]$ is a map which associates 
each simplicial complex $\Delta$ on $[n]$ with a simplicial complex 
$\Shift(\Delta)$ on $[n]$ 
and which satisfies the following conditions:
\begin{enumerate}
\item[(S$_1$)] $\Shift(\Delta)$ is shifted;
\item[(S$_2$)] $\Shift(\Delta) = \Delta$ if $\Delta$ is shifted;
\item[(S$_3$)] $f(\Delta) = f(\Shift(\Delta))$;
\item[(S$_4$)] $\Shift(\Delta') \subset \Shift(\Delta)$
if $\Delta' \subset \Delta$.
\end{enumerate}

\medskip

Erd\"os, Ko and Rado \cite{EKR} introduce a combinatorial shifting.
Let $\Delta$ be a simplicial complex on $[n]$.
Let $1 \leq i < j \leq n$.
Write
$\Shift_{ij}(\Delta)$
for the simplicial complex on $[n]$   
whose faces are 
$C_{ij}(\sigma) \subset [n]$,
where $\sigma \in \Delta$ and where
\begin{eqnarray*}
C_{ij}(\sigma) =
\left\{ 
\begin{array}{l}
(\sigma \setminus \{ i \}) \bigcup \{ j \}, 
\, \, \, \, \, \mbox{if} \, \, \, 
i \in \sigma, \, \, \,  
j \not\in \sigma \, \, \, 
\mbox{and} \, \, \, 
(\sigma \setminus \{ i \}) \bigcup \{ j \} \not\in \Delta,
\\
\hspace{0cm} \sigma, \hspace{2.6cm}
\mbox{otherwise.} 
\end{array} 
\right.
\end{eqnarray*} 
It follows from, e.g., \cite[Corollary 8.6]{Herzog}
that there exists a finite sequence of pairs of integers
$(i_1, j_1), (i_2, j_2), \ldots, (i_q, j_q)$
with each $1 \leq i_k < j_k \leq n$ such that
\[
\Shift_{i_q j_q} (\Shift_{i_{q-1} j_{q-1}}( \cdots 
(\Shift_{i_1 j_1}(\Delta)) \cdots ))
\]
is shifted.  Such a shifted complex is called 
a {\em combinatorial shifted complex} of $\Delta$ 
and will be denoted by $\Delta^c$. 
A combinatorial shifted complex $\Delta^c$ of $\Delta$ is, 
however, not necessarily unique.
The shifting operation $\Delta \mapsto \Delta^c$, which is 
a shifting operation (\cite[Lemma 8.4]{Herzog}),
is called {\em combinatorial shifting}.

Assume that the base field $K$ is infinite.
The {\em exterior algebraic shifted complex} of 
a simplicial complex $\Delta$ on $[n]$ is
the simplicial complex $\Delta^e$ on $[n]$ with
\[
J_{\Delta^e} = \Gin^E(J_\Delta).
\]
Following \cite[p.\ 105]{Herzog} and \cite[p.\ 125]{Kalai}
the shifting operation $\Delta \mapsto \Delta^e$, which is 
a shifting operation (\cite[Proposition 8.8]{Herzog}),
is called {\em exterior algebraic shifting}.

Assume that the base field $K$ is of
characteristic $0$.
Let $\Delta$ be a simplicial complex on $[n]$
and write $G(\Gin^S(I_\Delta))$ for 
the unique minimal system of monomial generators
of the generic initial ideal $\Gin^S(I_\Delta)$
of the Stanley--Reisner ideal $I_\Delta$ of $S$. 
Let $u = x_{i_1} x_{i_2} \cdots x_{i_j} \cdots x_{i_s}$,
where $1 \leq i_1 \leq i_2 \leq \cdots \leq i_j \leq 
\cdots \leq i_s \leq n$,
be a monomial belonging to $G(\Gin^S(I_\Delta))$.
One has $i_s + (s - 1) \leq n$
(\cite[Lemma 8.15]{Herzog}).
We then introduce the squarefree monomial
\[
u^* = x_{i_1} x_{i_2 + 1} \cdots x_{i_j + (j - 1)} \cdots 
x_{i_s + (s - 1)}
\]
of $S$ and write $(\Gin^S(I_\Delta))^*$ 
for the squarefree ideal of $S$
generated by those monomials $u^*$ with 
$u \in G(\Gin^S(I_\Delta))$.  
The {\em symmetric algebraic shifted complex}
of $\Delta$ is the simplicial complex $\Delta^s$ on 
$[n]$ with
\[
I_{\Delta^s} = (\Gin^S(I_\Delta))^*
\]
Since $\Gin^S(I_\Delta)$ is strongly stable,
it follows that $\Delta^s$ is shifted
(\cite[Lemma 8.17]{Herzog}).
The shifting operation $\Delta \mapsto \Delta^s$, 
which is a shifting operation 
(\cite[Theorem 8.19]{Herzog}),
is called {\em symmetric algebraic shifting}.

\section{Graded Betti numbers of
$I_{\Delta^e}$ and $I_{\Delta^c}$}
Let $K$ be an infinite field, 
$S = K[x_1, \ldots, x_n]$  
the polynomial ring in $n$ variables over $K$ with
each $\deg x_i = 1$
and $E = \bigoplus_{d=0}^{n} \bigwedge^d(V)$ 
the exterior algebra of a vector space $V$ over $K$
of dimension $n$ with basis $e_1, \ldots, e_n$. 
Assume that the general linear group
$\GL(n;K)$ acts linearly on $E$.
Let, as before, 
$<_{\rev}$ be the reverse lexicographic order on $E$ 
induced by the ordering $e_1 > \cdots > e_n$.

Given an arbitrary homogeneous ideal $I = \bigoplus_{d=0}^{n} I_d$
of $E$ with each $I_d \subset \bigwedge^d(V)$,
fix $\varphi \in \GL(n;K)$ for which 
$\initial_{<_{\rev}}(\varphi(I))$
is the generic initial ideal $\Gin^E(I)$ of $I$.
Recall that the subspace $\bigwedge^d(V)$
is of dimension ${n \choose d}$ 
with a canonical $K$-basis 
$e_\sigma$, $\sigma \in {[n] \choose d}$.
Choose an arbitrary $K$-basis $f_1, \ldots, f_s$ of  
$I_d$, where $s = \dim_K I_d$.  Write each $\varphi(f_i)$,
$1 \leq i \leq s$, of the form
\[
\varphi(f_i) = \sum_{\sigma \in {[n] \choose d}} 
\alpha_{i}^{\sigma} \, e_{\sigma}
\]
with each $\alpha_{i}^{\sigma} \in K$.
Let $M(I,d)$ denote the $s \times {n \choose d}$ matrix
\[
M(I,d) =(\alpha_{i}^{\sigma})_{1 \leq i \leq s, \,
\sigma \in {[n] \choose d}}
\]
whose columns are indexed by 
$\sigma \in {[n] \choose d}$.
Moreover, for each $\tau \in {[n] \choose d}$,
write 
$M_{\tau}(I,d)$ for the submatrix of $M(I,d)$
which consists of the columns of $M(I,d)$ indexed by
those $\sigma \in {[n] \choose d}$ with 
$e_\tau \leq_{\rev} e_\sigma$
and write
$M'_{\tau}(I,d)$ for the submatrix of $M_{\tau}(I,d)$
which is obtained by removing the column 
of $M_{\tau}(I,d)$ indexed by $\tau$.

\begin{Lemma}
\label{generic}
Let $e_{\tau} \in \bigwedge^d(V)$
with $\tau \in {[n] \choose d}$. 
Then
one has $e_{\tau} \in (\Gin^E(I))_d$
if and only if
$\rank(M'_{\tau}(I,d)) < \rank(M_{\tau}(I,d))$.
\end{Lemma}

\begin{pf}
In linear algebra
we know that
$\rank(M'_{\tau}(I,d)) < \rank(M_{\tau}(I,d))$
if and only if
the row vector $(0, \ldots, 0, 1)$
with ``$1$'' lying on the column indexed by $\tau$ 
arises in $M_{ \tau}(I,d)$
after repeating the elementary transformations 
on the row vectors of $M_{\tau}(I,d)$.
Thus, by identifying the rows of $M(I,d)$
with $\varphi(f_1), \ldots, \varphi(f_s)$,
it follows that 
$\rank(M'_{\tau}(I,d)) < \rank(M_{\tau}(I,d))$
if and only if there exist
$c_1, \ldots, c_{s}$ belonging to $K$
with
$\initial_{<_{\rev}}(f) = e_{\tau}$, where
$f = \sum_{i=1}^{s} c_i \varphi(f_i) \in (\varphi(I))_d$.
Since $\Gin^E(I) = \initial_{<_{\rev}}(\varphi(I))$,
one has $e_{\tau} \in (\Gin^E(I))_d$
if and only if
$\rank(M'_{\tau}(I,d)) < \rank(M_{\tau}(I,d))$,
as desired.
\end{pf}

\begin{Corollary}
\label{independent}
The rank of a matrix $M_{\tau}(I,d)$,
$\tau \in {[n] \choose d}$,
is independent of  
the choice of $\varphi \in \GL(n;K)$
for which 
$\Gin^E(I) = \initial_{<_{\rev}}(\varphi(I))$
together with
a $K$-basis $f_1, \ldots, f_s$ of $I_d$.
\end{Corollary}

\begin{Corollary}
\label{isomorphic}
Let $I \subset E$ be a homogeneous ideal
and $\psi \in \GL(n;K)$. 
Then one has 
$\rank(M_{\tau}(I,d)) = \rank(M_{\tau}(\psi(I),d))$
for all $\tau \in {[n] \choose d}$. 
\end{Corollary}


\begin{pf}
Recall that
there is a nonempty subset $U \subset \GL(n;K)$
which is Zariski open and dense such that
$\Gin^E(I) = \initial_{<_{\rev}}(\varphi(I))$
for all $\varphi \in U$.  Similarly,
there is a nonempty subset $V \subset \GL(n;K)$
which is Zariski open and dense such that
$\Gin^E(\psi(I)) = \initial_{<_{\rev}}(\varphi'(\psi(I)))$
for all $\varphi' \in V$.  
Since $U \psi^{-1} \bigcap V \neq \emptyset$, if
$\rho \in U \psi^{-1} \bigcap V$,
then $\Gin^E(I) = \initial_{<_{\rev}}(\rho(\psi(I))
= \Gin^E(\psi(I))$ and the matrix
$M(I,d)$ with using 
$\rho \psi \in U$ and a $K$-basis $f_1, \ldots, f_s$ of $I_d$
coincides with 
$M(\psi(I),d)$ with using $\rho \in V$
and a $K$-basis $\psi(f_1), \ldots, \psi(f_s)$ of $\psi(I)_d$.
\end{pf}

If $u = e_\sigma$ 
is a monomial of $E$, then we set 
$m(u) = \max\{ \, j \, : \, j \in \sigma \, \}$.
Given a monomial ideal $I \subset E$, 
one defines 
$m_{\leq i}(I, d)$,
where $1 \leq i \leq n$ and $1 \leq d \leq n$, by
\begin{eqnarray*}
m_{\leq i}(I, d) =
|\{ \, u = e_\sigma \in I \, : \, \deg(u) = d, m(u) \leq i \, \}|.
\end{eqnarray*}

\begin{Corollary}
\label{rank}
Let
$\sigma_{(i,d)} = \{ i - d + 1, i - d + 2, \ldots, i \}
\in {[n] \choose d}$.
Then given a homogeneous ideal
$I \subset E$ one has
\[
m_{\leq i}(\Gin^E(I), d) = \rank (M_{\sigma_{(i,d)}}(I,d)),
\]
where $\rank (M_{\sigma_{(i,d)}}(I,d)) = 0$
if $i < d$.
\end{Corollary}

\begin{pf}
Let $\tau \in {[n] \choose d}$.  Then
$m(e_\tau) \leq i$ if and only if
$e_{\sigma_{(i,d)}} \leq_{\rev} e_\tau$.
On the other hand,
Lemma \ref{generic} says that 
$\rank (M_{\sigma_{(i,d)}}(I,d))$
coincides with the number of monomials
$e_\tau \in (\Gin^E(I))_d$ with
$e_{\sigma_{(i,d)}} \leq_{\rev} e_\tau$.
Thus 
$m_{\leq i}(\Gin^E(I), d) = \rank (M_{\sigma_{(i,d)}}(I,d))$,
as required.
\end{pf}

Let $I \subset E$ be a monomial ideal.
Fix $1 \leq i < j \leq n$.
Let $t \in K$ and introduce the linear injective map 
$S_{ij}^t : I \to E$ satisfying
\begin{eqnarray*}
S_{ij}^t(e_\sigma) =
\left\{ 
\begin{array}{l}
e_{(\sigma \setminus \{ j \}) \bigcup \{ i \}}
+ t e_\sigma, 
\, \, \, \, \, \mbox{if} \, \, \, 
j \in \sigma, \, \, \,  
i \not\in \sigma \, \, \, 
\mbox{and} \, \, \, 
e_{(\sigma \setminus \{ j \}) \bigcup \{ i \}} \not\in I,
\\
\hspace{0cm} e_{\sigma}, \hspace{2.9cm}
\mbox{otherwise,} 
\end{array} 
\right.
\end{eqnarray*} 
where $e_\sigma \in I$ is a monomial.  
Let $I_{ij}(t) \subset E$ denote the image
of $I$ by $S_{ij}^t$.

\begin{Lemma}
\label{subspace}
{\em (a)}  
If $t \neq 0$, then there is 
$\lambda_{ij}^t \in GL(n;K)$ with
$I_{ij}(t) = \lambda_{ij}^t(I)$.
In particular 
the subspace $I_{ij}(t)$ is an ideal of $E$.

{\em (b)} 
Let $\Delta$ denote the simplicial complex
on $[n]$ and $J_\Delta$ its exterior face ideal.
Then 
$(J_\Delta)_{ij}(0) = J_{\Shift_{ij}(\Delta)}$.
\end{Lemma}

\begin{pf}
(a)  Let $\lambda_{ij}^t \in GL(n;K)$ satisfy
\begin{eqnarray*}
\lambda_{ij}^t(e_k) =
\left\{ 
\begin{array}{l} 
e_k 
\hspace{2cm}
(k \neq j),
\\
e_i + t e_j 
\hspace{1.05cm}
 (k = j)
\end{array} 
\right.
\end{eqnarray*}
We claim $I_{ij}(t) = \lambda_{ij}^t(I)$.
Let $e_\sigma \in I$.
\begin{enumerate}
\item[(i)]
If $j \not\in \sigma$, then
$\lambda_{ij}^t(e_\sigma) = e_\sigma
= S_{ij}^t(e_\sigma)$.
Thus $\lambda_{ij}^t(e_\sigma) \in I_{ij}(t)$.
\item[(ii)]
If
$j \in \sigma$ and $i \in \sigma$, then
$\lambda_{ij}^t(e_\sigma) = t e_\sigma
= t S_{ij}^t(e_\sigma)$.
Thus 
$\lambda_{ij}^t(e_\sigma) \in I_{ij}(t)$.
\item[(iii)]
Let
$j \in \sigma$ and $i \not\in \sigma$
with $e_{(\sigma \setminus \{ j \}) \bigcup \{ i \}}
\in I$.  Then 
$\lambda_{ij}^t(e_\sigma) = 
e_{(\sigma \setminus \{ j \}) \bigcup \{ i \}}
+ t e_\sigma$
and
$S_{ij}^t(e_\sigma) = e_\sigma$.
Since 
$e_{(\sigma \setminus \{ j \}) \bigcup \{ i \}}
\in I$,
$S_{ij}^t(e_{(\sigma \setminus \{ j \}) \bigcup \{ i \}})
= e_{(\sigma \setminus \{ j \}) \bigcup \{ i \}}
\in I_{ij}(t)$.
Thus
$\lambda_{ij}^t(e_\sigma) \in I_{ij}(t)$.
\item[(iv)]
Let
$j \in \sigma$ and $i \not\in \sigma$
with $e_{(\sigma \setminus \{ j \}) \bigcup \{ i \}}
\not\in I$.
Then
$\lambda_{ij}^t(e_\sigma) = 
e_{(\sigma \setminus \{ j \}) \bigcup \{ i \}}
+ t e_\sigma$
and
$S_{ij}^t(e_\sigma) = 
e_{(\sigma \setminus \{ j \}) \bigcup \{ i \}}
+ t e_\sigma$.
Thus
$\lambda_{ij}^t(e_\sigma) \in I_{ij}(t)$.
\end{enumerate}
Hence $\lambda_{ij}^t(I) \subset I_{ij}(t)$.
Since each of $\lambda_{ij}^t$ and $S_{ij}^t$
is injective, one has
$I_{ij}(t) = \lambda_{ij}^t(I)$,
as desired.

(b)  We claim 
$\{ \, \sigma \subset [n] \, : \, 
e_\sigma \in (J_\Delta)_{ij}(0) \, \}
\bigcap
\Shift_{ij}(\Delta)
= \emptyset$.
\begin{enumerate}
\item[(i)]
If $e_\sigma \in (J_\Delta)_{ij}(0)$
with $e_\sigma \not\in J_\Delta$, then 
there is $e_\tau \in J_\Delta$ 
with 
$\sigma = (\tau \setminus \{ j \}) \bigcup \{ i \}$.
Since $\sigma \in \Delta$, $\tau \not\in \Delta$
and
$\tau = (\sigma \setminus \{ i \}) \bigcup \{ j \}$,
one has 
$\tau = C_{ij}(\sigma) \in \Shift_{ij}(\Delta)$.
Thus $\sigma \not\in \Shift_{ij}(\Delta)$.
\item[(ii)]
Let $e_\sigma \in (J_\Delta)_{ij}(0)$
with $e_\sigma \in J_\Delta$.
Suppose $\sigma \in \Shift_{ij}(\Delta)$.
Since $\sigma \not\in \Delta$,
there is $\tau \subset [n]$ with
$\tau \in \Delta$ such that
$\sigma = (\tau \setminus \{ i \}) \bigcup \{ j \}$.
Hence
$j \in \sigma$, $i \not\in \sigma$ and
$e_\tau = e_{(\sigma \setminus \{ j \}) \bigcup \{ i \}}
\not\in J_\Delta$.
Thus $e_\tau \in (J_\Delta)_{ij}(0)$
and $e_\sigma \not\in (J_\Delta)_{ij}(0)$.
\end{enumerate}
Hence $(J_\Delta)_{ij}(0)  
\subset J_{\Shift_{ij}(\Delta)}$.
Since $\dim_{K} (J_\Delta)_{ij}(0)
= \dim_K J_\Delta
= \dim_{K} J_{\Shift_{ij}(\Delta)}$,
it follows that
$(J_\Delta)_{ij}(0)
= J_{\Shift_{ij}(\Delta)}$.
\end{pf}

\begin{Lemma}
\label{excellent}
Work with the same notation as in Corollary \ref{rank}.
One has
\[
\rank (M_{\sigma_{(i,d)}}(J_{\Shift_{ij}(\Delta)},d))
\leq
\rank (M_{\sigma_{(i,d)}}(J_\Delta,d)).
\]
\end{Lemma}

\begin{pf}
Fix a finite set $A \subset K$ with $0 \in A$ for which
$|A| \geq {[n] \choose d} + 2$.
One has $\varphi \in \GL(n;K)$
for which $\initial_{<_{\rev}}(\varphi((J_\Delta)_{ij}(t)))$
is the generic initial ideal of 
$(J_\Delta)_{ij}(t)$ for all $t \in A$.
For each $\sigma \in {[n] \choose d}$ we write 
\[
\varphi(e_\sigma) 
= \sum_{\tau \in {[n] \choose d}} c_\sigma^\tau e_\tau,
\, \, \, \, \, \, \, \, \, \, 
c_\sigma^\tau \in K.
\]
By using $\varphi$ together with the $K$-basis 
$\{ S_{ij}^t(e_\sigma) \, : \, e_\sigma \in (J_\Delta)_d \}$
of $((J_\Delta)_{ij}(t))_d$, we compute the matrix 
$M((J_\Delta)_{ij}(t),d)$.
If $S_{ij}^t(e_\sigma) =
e_{(\sigma \setminus \{ j \}) \bigcup \{ i \}} 
+ t e_\sigma$, then
\[
\varphi(S_{ij}^t(e_\sigma))
= \sum_{\tau \in {[n] \choose d}} 
(c^\tau_{(\sigma \setminus \{ j \}) \bigcup \{ i \}} 
+ t c^\tau_{\sigma}) e_\tau.
\]
Hence
\[
M((J_\Delta)_{ij}(t),d)
= (\alpha_\ell^\sigma + t \beta_\ell^\sigma)_{
1 \leq \ell \leq \dim_K((J_\Delta)_{ij}(t))_d, \,
\sigma \in {[n] \choose d}
}
\]
with each $\alpha_\ell^\sigma,
\beta_\ell^\sigma \in K$.

Let
$r(t) = 
\rank
(M_{\sigma_{(i,d)}}((J_\Delta)_{ij}(t),d))$. 
Thus $r(t)$ coincides with
the largest size of nonzero minors of the matrix
$M_{\sigma_{(i,d)}}((J_\Delta)_{ij}(t),d)$.
Fix a minor $N(t)$
of size $r(0)$ of
$M_{\sigma_{(i,d)}}((J_\Delta)_{ij}(t),d)$
with $N(0) \neq 0$.  We regard
$N(t)$ as a polynomial in $t$ of degree at most $r(0)$.
Since $r(0) \leq {[n] \choose d}$ and 
$|A| \geq {[n] \choose d} + 2$,
it follows that
there is $0 \neq a \in A$ with $N(a) \neq 0$.
Hence $r(0) \leq r(a)$.
Corollary \ref{isomorphic} 
together with Lemma \ref{subspace}
now guarantees that
$r(0) = \rank(M_{\sigma_{(i,d)}}
(J_{\Shift_{ij}(\Delta)},d))$
and
$r(a) = \rank (M_{\sigma_{(i,d)}}(J_\Delta,d))$.
Thus 
$\rank (M_{\sigma_{(i,d)}}(J_{\Shift_{ij}(\Delta)},d))
\leq
\rank (M_{\sigma_{(i,d)}}(J_\Delta,d))$,
as desired.
\end{pf}

\begin{Corollary}
\label{main}
Let $\Delta$ be a simplicial complex on $[n]$.
Then for all $i$ and $d$ one has
\[
m_{\leq i}(J_{\Delta^e},d) \geq m_{\leq i}(J_{\Delta^c},d).
\]
\end{Corollary}

\begin{pf}
Corollary \ref{rank} together with Lemma \ref{excellent}
guarantees that
\begin{eqnarray}
m_{\leq i}(\Gin^E(J_{\Delta}),d) \geq 
m_{\leq i}(\Gin^E(J_{\Shift_{ij}(\Delta)}),d).
\end{eqnarray}
Hence
$m_{\leq i}(\Gin^E(J_{\Delta}),d) \geq
m_{\leq i}(\Gin^E(J_{\Delta^c}),d)$.
In other words, one has
$m_{\leq i}(J_{\Delta^e},d) \geq
m_{\leq i}(J_{(\Delta^c)^e},d)$.
However, since $\Delta^c$ is shifted,
it follows that
$(\Delta^c)^e = \Delta^c$.
Thus
$m_{\leq i}(J_{\Delta^e},d) \geq
m_{\leq i}(J_{\Delta^c},d)$,
as desired.
\end{pf}

We now approach to the final step to prove 
the inequalities
$\beta_{ii+j}(I_{\Delta^e}) \leq \beta_{ii+j}(I_{\Delta^c})$
for all $i$ and $j$
on graded Betti numbers of $I_{\Delta^e}$ and $I_{\Delta^c}$.
Lemma \ref{formula} stated below 
essentially appears in 
\cite[pp.\ 376 -- 377]{AraHerHiLexsegment}.

\begin{Lemma}
\label{formula}
If $\Delta$ is a shifted simplicial complex, then
for all $i$ and $j$ one has
\[
\beta_{i i+j}(I_\Delta) = 
m_{\leq n}(I_\Delta,j) {{n-j} \choose {i}}
- \sum_{k=j}^{n-1} m_{\leq k}(I_\Delta,j) {{k-j} \choose {i-1}}
- \sum_{k=j}^{n} m_{\leq k-1}(I_\Delta,j-1) {{k-j} \choose i}.
\]
\end{Lemma}

\begin{Corollary}
\label{BostonMA}
Let $\Delta$ and $\Delta'$ be shifted simplicial complexes
on $[n]$ with $f(\Delta) = f(\Delta')$ and suppose that
\[
m_{\leq i}(J_\Delta, j) \geq m_{\leq i}(J_{\Delta'}, j)
\]
for all $i$ and $j$.  Then for all $i$ and $j$ one has
\[
\beta_{i i+j}(I_\Delta) 
\leq \beta_{i i+j}(I_{\Delta'}).
\]
\end{Corollary}

\begin{pf}
Since $f(\Delta) = f(\Delta')$, one has
$m_{\leq n}(I_\Delta, j) 
= m_{\leq n}(I_{\Delta'}, j)$
for all $j$.
Lemma \ref{formula} then yields 
the inequalities $\beta_{i i+j}(I_\Delta) 
\leq \beta_{i i+j}(I_{\Delta'})$
for all $i$ and $j$, as desired.
\end{pf}

\begin{Theorem}
\label{maintheorem}
Let the base field be infinite. 
Let $\Delta$ be a simplicial complex,
$\Delta^e$ the exterior algebraic shifted complex of
$\Delta$ and
$\Delta^c$ a combinatorial shifted complex of
$\Delta$.
Then 
\[
\beta_{ii+j}(I_{\Delta^e}) \leq \beta_{ii+j}(I_{\Delta^c})
\] 
for all $i$ and $j$. 
\end{Theorem}

\begin{pf}
Corollary \ref{main} guarantees 
$m_{\leq i}(J_{\Delta^c}, j) \leq m_{\leq i}(J_{\Delta^e}, j)$
for all $i$ and $j$.
Thus by virtue of Corollary \ref{BostonMA}
the required inequalities
$\beta_{ii+j}(I_{\Delta^e}) \leq \beta_{ii+j}(I_{\Delta^c})$
follow immediately.
\end{pf}

\section{Graded Betti numbers of
$I_{\Delta}$ and $I_{\Delta^c}$}
Let $K$ be an arbitrary field, and let
$S = K[x_1, \ldots, x_n]$ denote the polynomial ring 
in $n$ variables
over $K$ with each $\deg x_i = 1$.
Let $\Delta$ be a simplicial complex on $[n]$ and
$I_\Delta \subset S$ its Stanley--Reisner ideal.
Let ${\tilde H}_k(\Delta;K)$ denote the $k$th reduced 
homology group of $\Delta$ with coefficients $K$.
If $W \subset [n]$, then $\Delta_W$ stands for 
the simplicial complex on $W$ 
whose faces are those faces $\sigma$ of $\Delta$ 
with $\sigma \subset W$.

Recall that Hochster's formula 
\cite[Theorem 5.5.1]{BrunsHerzog}
to compute the graded Betti numbers of $I_\Delta$ says that
\begin{eqnarray}
\beta_{ii+j}(I_\Delta)
= \sum_{W \subset [n], \, |W|=i+j} 
\dim_K({\tilde H}_{j-2}(\Delta_W;K))
\end{eqnarray}
for all $i$ and $j$.

Fix $1 \leq i < j \leq n$ and set 
$\Gamma = \Shift_{ij}(\Delta)$.

\begin{Lemma}
\label{Boston}
One has
\[
\dim_K({\tilde H}_{k}(\Delta;K))
\leq 
\dim_K({\tilde H}_{k}(\Gamma;K))
\]
for all $k$.
\end{Lemma}

\begin{pf}
By considering an extension field of $K$ if necessarily,
we assume that $K$ is infinite. 
Let $\Delta^e$ denote 
the exterior algebraic shifted complex  
of $\Delta$.  
It is known \cite[Proposition 8.10]{Herzog}
that  
${\tilde H}_{k}(\Delta;K) \cong {\tilde H}_{k}(\Delta^e;K)$.
Thus what we must prove is 
$\dim_K({\tilde H}_{k}(\Delta^e;K))
\leq 
\dim_K({\tilde H}_{k}(\Gamma^e;K))$
for all $k$. 
By using $(2)$ one has
$\beta_{i n}(I_\Delta) 
= \dim_K({\tilde H}_{n-i-2}(\Delta;K))$.
Hence our work is to show that
$\beta_{i n}(I_{\Delta^e}) \leq \beta_{i n}(I_{\Gamma^e})$
for all $i$.
The inequality $(1)$ says that
$m_{\leq i}(J_{\Delta^e},j) \geq m_{\leq i}(J_{\Gamma^e},j)$
for all $i$ and $j$.
It then follows from Corollary \ref{BostonMA} that
$\beta_{i i+j}(I_{\Delta^e}) 
\leq \beta_{i i+j}(I_{\Gamma^e})$
for all $i$ and $j$.
Thus in particular
$\beta_{i n}(I_{\Delta^e}) \leq \beta_{i n}(I_{\Gamma^e})$
for all $i$.
\end{pf}

Let $W \subset [n] \setminus \{ i, j \}$.
Let $\Delta_1 = \Delta_{W \cup \{ i \}}$,
$\Delta_2 = \Delta_{W \cup \{ j \}}$,
$\Gamma_1 = \Gamma_{W \cup \{ i \}}$
and
$\Gamma_2 = \Gamma_{W \cup \{ j \}}$.
Then 
\begin{eqnarray*}
\Delta_1 \bigcap \Delta_2 = \Gamma_1 \bigcap \Gamma_2 
= \Delta_W = \Gamma_W,
\end{eqnarray*}
\begin{eqnarray}
\Gamma_1 \bigcup \Gamma_2 = \Shift_{ij}(\Delta_1 \bigcup \Delta_2).
\end{eqnarray}
Recall that the reduced Mayer--Vietoris exact sequence 
of $\Delta_1$ and $\Delta_2$ 
and that of $\Gamma_1$ and $\Gamma_2$ are
the exact sequences
\begin{eqnarray*}
\def\normalbaseline{\baselineskip20pt
     \lineskip3pt  \lineskiplimit3pt}
\def\mapright#1{\smash{
     \mathop{\longrightarrow}\limits^{#1}}}
\cdots
\mapright{}  
{\tilde H}_k(\Delta_W;K)
\mapright{\partial_{1,k}}
{\tilde H}_k(\Delta_1;K) \bigoplus {\tilde H}_k(\Delta_2;K)
\mapright{\partial_{2,k}}
{\tilde H}_k(\Delta_1 \bigcup \Delta_2;K)
\\  
\def\normalbaseline{\baselineskip20pt
     \lineskip3pt  \lineskiplimit3pt}
\def\mapright#1{\smash{
     \mathop{\longrightarrow}\limits^{#1}}}
\mapright{\partial_{3,k}}  
{\tilde H}_{k-1}(\Delta_W;K)
\mapright{\partial_{1,k-1}}
\cdots,
\hspace{6.75cm}
\end{eqnarray*}
\begin{eqnarray*}
\def\normalbaseline{\baselineskip20pt
     \lineskip3pt  \lineskiplimit3pt}
\def\mapright#1{\smash{
     \mathop{\longrightarrow}\limits^{#1}}}
\cdots
\mapright{}
{\tilde H}_k(\Gamma_W;K)
\mapright{\partial'_{1,k}}
{\tilde H}_k(\Gamma_1;K) \bigoplus {\tilde H}_k(\Gamma_2;K)
\mapright{\partial'_{2,k}}
{\tilde H}_k(\Gamma_1 \bigcup \Gamma_2;K) 
\\
\def\normalbaseline{\baselineskip20pt
     \lineskip3pt  \lineskiplimit3pt}
\def\mapright#1{\smash{
     \mathop{\longrightarrow}\limits^{#1}}}
\mapright{\partial'_{3,k}}
{\tilde H}_{k-1}(\Gamma_W;K)
\mapright{\partial'_{1,k-1}}
\cdots.
\hspace{6.45cm}
\end{eqnarray*}

\begin{Lemma}
\label{kernel}
One has
\[
\Ker(\partial'_{1,k}) \subset \Ker(\partial_{1,k}).
\]
for all $k$.
\end{Lemma}

\begin{pf}
Let $\pi$ be a permutation on $[n]$ 
with $\pi(i) < \pi(j)$ and 
$\pi(\Delta)$ the simplicial complex
$\{ \pi(F) \, : \, F \in \Delta \}$ on $[n]$.
Since the combinatorial type of $\Shift_{ij}(\Delta)$
is equal to that of $\Shift_{\pi(i)\pi(j)}(\pi(\Delta))$,
we will assume that $j = i + 1$.

Let, in general, 
$C_k(\Delta)$ denote the vector space over $K$
with basis $\{ e_{i_0 i_1 \cdots i_k} \}$,
where 
$\{ i_0, i_1, \ldots, i_k \} \in \Delta$
and where
$1 \leq i_0 < i_1 < \cdots < i_k \leq n$,
and define the linear map
$\partial : C_k(\Delta) \to C_{k-1}(\Delta)$
by setting $\partial(e_{i_0 i_1 \cdots i_k})
= \sum_{j=0}^{k} 
( - 1 )^j e_{i_0 i_1 \cdots i_{j-1} i_{j+1} \cdots i_k}$.

Let $[a] \in \Ker(\partial'_{1,k})$, where
$a \in C_k(\Gamma_W)$.  Since 
$([a], [a])  
\in {\tilde H}_k(\Gamma_1;K) \bigoplus {\tilde H}_k(\Gamma_2;K)$
vanishes, one has $u \in C_{k+1}(\Gamma_1)$ with
$\partial(u) = a$.  Say,
\[
u = \sum_{|F| = k + 1, \, i \not\in F, \, F \cup \{ i \} \in \Gamma_1} 
a_{F \cup \{ i \}} e_{F \cup \{ i \}} 
+ \sum_{|G| = k + 2, \, 
G \in \Delta_W} 
b_G e_G,
\]
where $a_{F \cup \{ i \}}, b_G \in K$. 

Let $F \subset W$ with $F \cup \{ i \} \in \Gamma_1$.
Then
$F \cup \{ i \} \in \Delta_1$
and 
$F \cup \{ j \} \in \Delta_2$.
Thus $F \cup \{ j \} \in \Gamma_2$. 
In particular $u \in C_{k+1}(\Delta_1)$
with $\partial(u) = a$.

Since $a \in C_k(\Gamma_W)$ 
is a linear combination of those basis elements $e_F$
with $F \in \Gamma$, $F \subset W$ and $|F| = k + 1$
and since $j = i + 1$,
it follows that $\partial(v) = a$, 
where $v \in C_{k+1}(\Gamma_2) \bigcap C_{k+1}(\Delta_2)$ 
is the element
\[
v = \sum_{|F| = k + 1, \, i \not\in F, \, F \cup \{ i \} \in \Gamma_1} 
a_{F \cup \{ i \}} e_{F \cup \{ j \}} 
+ \sum_{|G| = k + 2, \, 
G \in \Delta_W}  
b_G e_G.
\] 
Hence 
$([a], [a])  
\in {\tilde H}_k(\Delta_1;K) \bigoplus {\tilde H}_k(\Delta_2;K)$
vanishes, as required.
\end{pf}

It then follows that
\begin{eqnarray*}
\dim_K(\Ker(\partial_{1,k}))
\geq
\dim_K(\Ker(\partial'_{1,k})),
\end{eqnarray*}
\begin{eqnarray*}
\dim_K(\Im(\partial_{1,k}))
\leq
\dim_K(\Im(\partial'_{1,k})),
\end{eqnarray*}
\begin{eqnarray}
\dim_K(\Ker(\partial_{2,k}))
\leq
\dim_K(\Ker(\partial'_{2,k})).
\end{eqnarray}

On the other hand,  
\begin{eqnarray}
\dim_K({\tilde H}_k(\Delta_1 \bigcup \Delta_2;K)) 
= \dim_K(\Ker(\partial_{3,k})) + \dim_K(\Im(\partial_{3,k})),
\end{eqnarray}
\begin{eqnarray}
\dim_K({\tilde H}_k(\Gamma_1 \bigcup \Gamma_2;K)) 
= \dim_K(\Ker(\partial'_{3,k})) + \dim_K(\Im(\partial'_{3,k})).
\end{eqnarray}
Lemma \ref{Boston} together with $(3)$ guarantees that
\begin{eqnarray}
\dim_K({\tilde H}_k(\Delta_1 \bigcup \Delta_2;K)) 
\leq
\dim_K({\tilde H}_k(\Gamma_1 \bigcup \Gamma_2;K)). 
\end{eqnarray}
Since 
$\Im(\partial_{3,k}) = \Ker(\partial_{1,k-1})$
and
$\Im(\partial'_{3,k}) = \Ker(\partial'_{1,k-1})$,
Lemma \ref{kernel} yields 
\begin{eqnarray}
\dim_K(\Im(\partial_{3,k})) 
\geq
\dim_K(\Im(\partial'_{3,k})). 
\end{eqnarray}
Since 
$\Im(\partial_{2,k}) = \Ker(\partial_{3,k})$
and
$\Im(\partial'_{2,k}) = \Ker(\partial'_{3,k})$,
it follows from
$(5)$ and $(6)$ together with $(7)$ and $(8)$ 
that
\begin{eqnarray}
\dim_K(\Im(\partial_{2,k}))
\leq
\dim_K(\Im(\partial'_{2,k})). 
\end{eqnarray}

Finally, it follows from 
the reduced Mayer--Vietoris exact sequence 
of $\Delta_1$ and $\Delta_2$ 
and that of $\Gamma_1$ and $\Gamma_2$
together with $(4)$ and $(9)$ that 
\begin{eqnarray}
\dim_K({\tilde H}_{k}(\Delta_{1};K)
\bigoplus {\tilde H}_{k}(\Delta_{2};K)) 
\leq 
\dim_K({\tilde H}_{k}(\Gamma_{1};K)
\bigoplus {\tilde H}_{k}(\Gamma_{2};K)).
\end{eqnarray}

\begin{Lemma}
\label{Sydney}
Fix $1 \leq p < q \leq n$.
Let $\Delta$ be a simplicial complex on $[n]$
and $\Gamma = \Shift_{pq}(\Delta)$.
Then
\[
\beta_{i i+j}(I_\Delta) \leq \beta_{i i+j}(I_\Gamma)
\]
for all $i$ and $j$.
\end{Lemma}

\begin{pf}
The right-hand side of Hochster's formula $(11)$ 
can be rewritten as
\[
\beta_{i i+j}(I_\Delta) = \alpha_{ij}(\Delta) +
\gamma_{ij}(\Delta) + \delta_{ij}(\Delta),
\]
where
\begin{eqnarray*}
\alpha_{ij}(\Delta) 
& = & 
\sum_{W \subset [n] \setminus \{ p, q \}, \, |W|=i+j}
\dim_K({\tilde H}_{j-2}(\Delta_W;K)), \\
\gamma_{ij}(\Delta) 
& = &
\sum_{W \subset [n] \setminus \{ p, q \}, \, |W|=i+j-1} 
\dim_K({\tilde H}_{j-2}(\Delta_{W \cup \{ p \}};K)) \\
&   & \, \, \, \, \, \, \, \, \, \, +
\sum_{W \subset [n] \setminus \{ p, q \}, \, |W|=i+j-1} 
\dim_K({\tilde H}_{j-2}(\Delta_{W \cup \{ q \}};K)), \\
\delta_{ij}(\Delta)
& = & 
\sum_{W \subset [n] \setminus \{ p, q \}, \, |W|=i+j-2} 
\dim_K({\tilde H}_{j-2}(\Delta_{W \cup \{ p, q \}};K)).
\end{eqnarray*}

Let $W \subset [n] \setminus \{ p, q \}$.
Then $\Delta_W = \Gamma_W$.  Thus
$\alpha_{ij}(\Delta) = \alpha_{ij}(\Gamma)$.
Since
$\Gamma_{W \cup \{ p, q \}} =
\Shift(\Delta_{W \cup \{ p, q \}})$,
Lemma \ref{Boston} says that
$\delta_{ij}(\Delta) \leq \delta_{ij}(\Gamma)$.
Finally, it follows from $(10)$ that
$\gamma_{ij}(\Delta) \leq \gamma_{ij}(\Gamma)$.
Hence 
$\beta_{i i+j}(I_\Delta) \leq \beta_{i i+j}(I_\Gamma)$,
as desired.
\end{pf}

Lemma \ref{Sydney} together with the definition of
combinatorial shifting now guarantees that

\begin{Theorem}
\label{theorem}
Let the base field be arbitrary.
Let $\Delta$ be a simplicial complex 
and $\Delta^c$ a combinatorial shifted complex
of $\Delta$.  Then
\[
\beta_{i i+j}(I_\Delta) \leq \beta_{i i+j}(I_{\Delta^c})
\]
for all $i$ and $j$.
\end{Theorem}

Let $\Delta'$ be a shifted simplicial complex
with the same $f$-vector as $\Delta$ and $\Delta^{\lex}$ 
the unique lexsegment simplicial complex
with the same $f$-vector as $\Delta$.
It is known \cite[Theorem 4.4]{AraHerHiLexsegment}
that
$\beta_{ii+j}(I_{\Delta'}) \leq
\beta_{ii+j}(I_{\Delta^{\lex}})$
for all $i$ and $j$.
Since $\Delta^c$ is shifted
with $f(\Delta^c) = f(\Delta)$,
it follows that
$\beta_{ii+j}(I_{\Delta^c}) \leq
\beta_{ii+j}(I_{\Delta^{\lex}})$
for all $i$ and $j$.
Hence 

\begin{Corollary}
Let the base field be arbitrary.
Let $\Delta$ be a simplicial complex 
and $\Delta^{\lex}$ 
the unique lexsegment simplicial complex
with the same $f$-vector as $\Delta$.
Then
\label{corollary}
\[
\beta_{ii+j}(I_\Delta) \leq
\beta_{ii+j}(I_{\Delta^{\lex}})
\]
for all $i$ and $j$.
\end{Corollary}

\section{Bad behavior of combinatorial shifted complexes}
Given a simplicial complex $\Delta$,
do there exist combinatorial shifted complexes
$\Delta_{\flat}^c$ and $\Delta_{\sharp}^c$ of $\Delta$ 
such that, for all combinatorial shifted complex
$\Delta^c$ of $\Delta$ and for all $i$ and $j$, one has
\[
\beta_{ii+j}(I_{\Delta_{\flat}^c}) 
\leq \beta_{ii+j}(I_{\Delta^c})
\leq \beta_{ii+j}(I_{\Delta_{\sharp}^c}) \, ?
\]
Unfortunately, in general,
the existence of such the combinatorial shifted complexes
$\Delta_{\flat}^c$ and $\Delta_{\sharp}^c$
cannot be expected.

Let $V$ be a vector space of dimension $15$ with basis
$e_1, \ldots, e_{15}$ 
and $E = \bigoplus_{d=0}^{15} \wedge^d(V)$
the exterior algebra of $V$. 
Let $<_{\lex}$ denote the lexicographic order on $E$
induced by the ordering 
$e_1 > \cdots > e_{15}$.
To simplify the notation we employ the following
\begin{eqnarray*}
h_1 = e_1, \, \, \, 
h_2 = e_2 \wedge e_3, \, \, \,
h_3 = e_3 \wedge e_4 \wedge e_5, 
\, \, \, \, \,  
\, \, \, \, \, \, \, \, \, \, \\
h_4 = e_4 \wedge \cdots \wedge e_{7}, \, \, 
h_5 = e_5 \wedge \cdots \wedge e_{9}, \, \, 
h_6 = e_6 \wedge \cdots \wedge e_{11}.
\end{eqnarray*}
First, we introduce $H_{i} \subset \wedge^2(V)$ 
with $3 \leq i \leq 8$ and $A, B \subset \wedge^2(V)$ 
by setting
\begin{eqnarray*}
H_{3} = \{ e_{12} \wedge e_{13}, e_{12} \wedge e_{15}, 
e_{13} \wedge e_{14} \}, & &
H_{4} = \{ e_{12} \wedge e_{13}, e_{12} \wedge e_{14}, 
e_{14} \wedge e_{15} \}, \\
H_{5} = \{ e_{12} \wedge e_{13}, e_{12} \wedge e_{15}, 
e_{14} \wedge e_{15} \}, & &
H_{6} = \{ e_{12} \wedge e_{13}, e_{13} \wedge e_{14}, 
e_{14} \wedge e_{15} \}, \\
H_{7} = \{ e_{12} \wedge e_{13}, e_{13} \wedge e_{15}, 
e_{14} \wedge e_{15} \}, & &
H_{8} = \{ e_{12} \wedge e_{14}, e_{13} \wedge e_{15}, 
e_{14} \wedge e_{15} \}, \\
A = \{ e_{12} \wedge e_{13} , e_{12} \wedge e_{14} , 
e_{13} \wedge e_{14} \}, & & \, \, \, 
B = \{ e_{12} \wedge e_{13} , e_{12} \wedge e_{14} , 
e_{12} \wedge e_{15} \}.
\end{eqnarray*}
Second, we introduce $T_i \subset \wedge^i(V)$
and $T_i(H) \subset \wedge^i(V)$
with $3 \leq i \leq 8$ by setting
\begin{eqnarray*}
T_i & = & \{ e_\sigma \in \wedge^i(V) \, : \, 
h_{i-2} \wedge e_{12} \wedge e_{13} <_{\lex} e_\sigma \}, \\ 
T_i(H) & = & \{ h_{i-2} \wedge e_\sigma \, : \, 
e_\sigma \in H \} \, \, \, \, \, \, \, \, \, \, 
\mbox{where} \, \, \, \, \, 
H \in \{ H_i, A, B \}.
\end{eqnarray*}
Let $I = \bigoplus_{d=3}^{15} I_d \subset E$ 
denote the ideal of $E$ 
generated by the monomials belonging to 
$\bigcup_{i=3}^{8} (T_i \bigcup T_i(H_i))$
together with all monomials of degree $9$
and $\Delta$ the simplicial complex on $\{ 1, \ldots, 15 \}$
with $I = J_\Delta$.

\begin{Lemma}
\label{routine}
{\em (a)}
For $3 \leq d \leq 8$ the subspace $I_d$ is spanned by
$T_d \bigcup T_d(H_d)$.

{\em (b)} 
Let $3 \leq d \leq 8$ and
$e_\sigma \in I_d$ with 
$e_\sigma \not\in T_d(H_d)$.
Then $S_{ij}^0(e_\sigma) = e_\sigma$.

{\em (c)} 
Unless $12 \leq i < j \leq15$
one has $S_{ij}^0(e_\sigma) = e_\sigma$
for all $e_\sigma \in \bigcup_{d=3}^{8} T_d(H_d)$.
\end{Lemma}

\begin{pf}
(a)  Let $3 \leq d < 8$.
We claim $e_j (T_d \bigcup T_d(H_d)) \subset T_{d+1}$
for all $j$.
In fact, $h_{d-1} \wedge e_{12} \wedge e_{13}
<_{\lex} 
e_j \wedge h_{d-2} \wedge e_{p} \wedge e_{q}$
unless 
$e_j \wedge h_{d-2} \wedge e_{p} \wedge e_{q} \neq 0$.

(b)  Let $e_\sigma \in I_d$
with $e_\sigma \not\in T_d(H_d)$. 
Let $j \in \sigma$ and $i \not\in \sigma$.
Since
$h_{d-2} \wedge e_{12} \wedge e_{13}
<_{\lex} e_\sigma$,
one has  
$h_{d-2} \wedge e_{12} \wedge e_{13}
<_{\lex} 
e_{(\sigma \setminus \{ j \}) \bigcup \{ i \}}$.
Thus 
$e_{(\sigma \setminus \{ j \}) \bigcup \{ i \}} 
\in T_d$.  Hence
$S_{ij}^0(e_\sigma) = e_\sigma$.

(c)  Let $i < 12$.
Let $e_\tau = h_{d-2} \wedge e_\sigma
\in T_d(H_d)$.
Let $j \in \tau$ and $i \not\in \tau$.
Then $h_{d-2} \wedge e_{12} \wedge e_{13}
<_{\lex} 
e_{(\tau \setminus \{ j \}) \bigcup \{ i \}}$.
Thus 
$e_{(\tau \setminus \{ j \}) \bigcup \{ i \}} \in T_d$.
Hence $S_{ij}^0(e_\sigma) = e_\sigma$.
\end{pf}

Given a sequence ${\bold Q} = (Q_3, \ldots, Q_8)$ with each
$Q_i \in \{ A, B \}$ we write $I^{\bold Q}$ for the ideal
of $E$ generated by the monomials belonging to 
$\bigcup_{i=3}^{8} (T_i \bigcup T_i(Q_i))$
together with all monomials of degree $9$.
Let ${\cal W}_{\shift}(\Delta)$ denote the set of shifted 
simplicial complexes of $\Delta$.

\begin{Lemma}
\label{computation}
{\em (a)}
Let $\Delta^c \in {\cal W}_{\shift}(\Delta)$.
Then $J_{\Delta^c}$ is of the form $I^{\bold Q}$.

{\em (b)}
None of $\Delta^c \in {\cal W}_{\shift}(\Delta)$ 
satisfies $J_{\Delta^c} = I^{(A, \ldots, A)}$.

{\em (c)}
None of $\Delta^c \in {\cal W}_{\shift}(\Delta)$ 
satisfies $J_{\Delta^c} = I^{(B, \ldots, B)}$.

{\em (d)}
For each $i$ and for each $j$ with $i < j$
there is $\Delta^c(i,j;A) \in {\cal W}_{\shift}(\Delta)$ 
with $J_{\Delta^c(i,j;A)} = I^{\bold Q}$, where
$Q_i = Q_j = A$.

{\em (e)}
For each $i$ and for each $j$ with $i < j$
there is $\Delta^c(i,j;B) \in {\cal W}_{\shift}(\Delta)$ 
with $J_{\Delta^c(i,j;B)} = I^{\bold Q}$, where
$Q_i = Q_j = B$.
\end{Lemma}

\begin{pf} 
After repeated applications of the operations
$S^0_{i_k j_k}$,
where $12 \leq i_k < j_k \leq 15$
and where $k = 1, 2, \ldots$,
each subset $T_d(H_d)$
will shift to 
either $T_d(A)$ or $T_d(B)$.
Moreover, 
$S_{ij}^0(T_d(A)) = T_d(A)$
and $S_{ij}^0(T_d(B)) = T_d(B)$
for all $1 \leq i < j \leq 15$.
Our claim (a) follows from
this observation together with
Lemma \ref{routine}.

A routine computation yields the
classification of the sequences
${\bold Q} = (Q_3, \ldots, Q_8)$ for which
there is 
$\Delta^c \in {\cal W}_{\shift}(\Delta)$
with
$J_{\Delta^c} = I^{\bold Q}$.
The classification table is 
\begin{eqnarray*}
(A,A,A,A,A,B), & & (A,A,A,A,B,A), \ldots, 
(B,A,A,A,A,A), \\
(B,B,B,B,B,A), & & (B,B,B,B,A,B), \ldots,  
(A,B,B,B,B,B)
\end{eqnarray*}
together with 
\begin{eqnarray*}
(A,A,A,B,B,B), & & (B,A,B,A,B,A), \\ 
(B,B,A,B,A,A), & & (A,B,B,A,A,B). 
\end{eqnarray*}
Our claims (b), (c), (d) and (e)
now follows immediately.
\end{pf}

\begin{Theorem}
\label{nonexistence}
{\em (a)}
None of $\Delta_{\sharp}^c \in {\cal W}_{\shift}(\Delta)$ satisfies 
$\beta_{ii+j}(J_{\Delta^c}) \leq \beta_{ii+j}(J_{\Delta_{\sharp}^c})$
for all $\Delta^c \in {\cal W}_{\shift}(\Delta)$
and for all $i$ and $j$.

{\em (b)}
None of $\Delta_{\flat}^c \in {\cal W}_{\shift}(\Delta)$ satisfies 
$\beta_{ii+j}(J_{\Delta_{\flat}^c}) \leq \beta_{ii+j}(J_{\Delta^c})$
for all $\Delta^c \in {\cal W}_{\shift}(\Delta)$
and for all $i$ and $j$.
\end{Theorem}

\begin{pf}
Let $\Delta_{\sharp}^c \in {\cal W}_{\shift}(\Delta)$
with $J_{\Delta_{\sharp}^c} = I^{{\bold Q}}$.
By Lemma \ref{computation} (c) 
there is $3 \leq j \leq 8$ with
$Q_j = A$ and $Q_{j'} = B$
for all $3 \leq j' < j$.
Lemma \ref{computation} (e) 
guarantees the existence of 
$\Delta^c(j-1,j;B) \in {\cal W}_{\shift}(\Delta)$ 
with $J_{\Delta^c(j-1,j;B)} = I^{{\bold Q}'}$, where
${\bold Q}' = (Q'_3, \ldots, Q'_8)$ with
$Q'_{j-1} = Q'_j = B$.
Then for $i \neq 14$ one has 
$m_{\leq i}(J_{\Delta^c(j-1,j;B)},j-1)
= m_{\leq i}(J_{\Delta_{\sharp}^c},j-1)$
and 
$m_{\leq i}(J_{\Delta^c(j-1,j;B)},j)
= m_{\leq i}(J_{\Delta_{\sharp}^c},j)$.
On the other hand, 
$m_{\leq 14}(J_{\Delta^c(j-1,j;B)},j-1)
= m_{\leq 14}(J_{\Delta_{\sharp}^c},j-1)$
and
$m_{\leq 14}(J_{\Delta^c(j-1,j;B)},j)
< m_{\leq 14}(J_{\Delta_{\sharp}^c},j)$.
Now, Lemma \ref{formula} says that
$\beta_{ii+j}(J_{\Delta_{\sharp}^c}) < 
\beta_{ii+j}(J_{\Delta^c(j-1,j;B)})$
for all $i$.
Thus $\Delta_\sharp^c \in {\cal W}_{\shift}(\Delta)$,
such that 
$\beta_{ii+j}(J_{\Delta^c}) 
\leq \beta_{ii+j}(J_{\Delta_{\sharp}^c})$
for all $\Delta^c \in {\cal W}_{\shift}(\Delta)$
and for all $i$ and $j$,
does not exist.
This completes the proof of (a).
Similar technique can be used to prove (b).
\end{pf}

\begin{Corollary}
\label{exteriorcombinatorial}
None of $\Delta^c \in {\cal W}_{\shift}(\Delta)$ satisfies
$\Delta^e = \Delta^c$. 
\end{Corollary}

\begin{pf}
Let
$\Delta_{\flat}^c \in {\cal W}_{\shift}(\Delta)$
satisfy $\Delta^e = \Delta_{\flat}^c$.
Since $\beta_{ii+j}(J_{\Delta^e}) \leq
\beta_{ii+j}(J_{\Delta^c})$
for all $i$ and $j$, it follows that
$\beta_{ii+j}(J_{\Delta_{\flat}^c}) \leq \beta_{ii+j}(J_{\Delta^c})$
for all $\Delta^c \in {\cal W}_{\shift}(\Delta)$
and for all $i$ and $j$.
This fact contradicts Theorem \ref{nonexistence} (b).
Thus none of $\Delta^c \in {\cal W}_{\shift}(\Delta)$ satisfies
$\Delta^e = \Delta^c$, as desired.
\end{pf}

\bigskip
\bigskip

{\small
\noindent
Satoshi Murai \\
Department of Pure and Applied Mathematics\\
Graduate School of Information Science and Technology\\
Osaka University \\
Toyonaka, Osaka 560-0043, Japan\\
E-mail:s-murai@@ist.osaka-u.ac.jp

\bigskip

\noindent
Takayuki Hibi\\
Department of Pure and Applied Mathematics\\
Graduate School of Information Science and Technology\\
Osaka University \\
Toyonaka, Osaka 560-0043, Japan\\
E-mail:hibi@@math.sci.osaka-u.ac.jp

}

\end{document}